Opinion

# Chance Constrained Optimization for Energy Management in Electric Vehicles

**Erfan Mohagheghi[1]*, Joan Gubianes Gasso[1], Abebe Geletu[2] and Pu Li[2]**

[1]MicroFuzzy GmbH, Taunusstraße 38, 80807 Munich, Germany

[2]Department of Process Optimization, Ilmenau University of Technology, Ilmenau, Germany

***Corresponding author:** Dr. Erfan Mohagheghi, MicroFuzzy GmbH, Taunusstraße 38, 80807 Munich, Germany, E-mail: Erfan.Mohagheghi@microfuzzy.com

E-powertrain of future electric vehicles could consist of energy generation units (e.g., fuel cells and photovoltaic modules), energy storage systems (e.g., batteries and supercapacitors), energy conversion units (e.g., bidirectional DC/DC converters and DC/AC inverters) and an electric machine, which can work in both generating and motoring modes [1-6]. An energy management system is responsible to operate the above-mentioned components in a way that the technical constraints are satisfied. This task should be accomplished by solving an optimization problem, which could aim at minimizing the total operation costs [5]. The optimization problem has been widely addressed by deterministic approaches [7], which take into account the forecasted values of active-reactive load profile. However, as shown in Figure 1 (a), it is impossible to accurately forecast the values, meaning that the solutions coming from deterministic approaches could lead to infeasible operations (i.e., constraint violations). Therefore, stochastic optimization approaches [8] should be utilized to find optimal solution strategies while considering uncertain parameters.

There are several mathematical approaches for optimization under uncertainty each of which could be suitable for a specific type of application. For instance, robust optimization and worst-case optimization is frequently used in many applications in which constraint violations are not tolerated [9]. However, in electric vehicles, there exist some types of constraints (e.g., upper and lower limit of state of the charge in batteries [10]), which are allowed to be violated to some degree and also for a limited time. For this, the method of chance constrained programming [11-13] (see Figure 1 (b)) can be used, by which optimal operations are obtained, and the satisfaction of constraints is ensured with a 'predefined probability' level. The deterministic equivalent of chance constraints can be determined easily if the model is linear and the uncertain values are 'normally' distributed. However, the e-powertrain model is mixed-integer nonlinear [14] and the uncertain variables are described by non-Gaussian probability density functions. In addition, differential equations of energy storage units introduce further complexities to the problem as stochastic 'dynamic' optimization problem has to be solved. For this, MicroFuzzy GmbH in collaboration with Ilmenau University of Technology formulates and solves this complex problem using the powerful methods [15,16] which have been recently developed in the Department of Process Optimization of that university. The solutions obtained by utilizing the chance constrained method could lead to significant reduction in total operation costs.

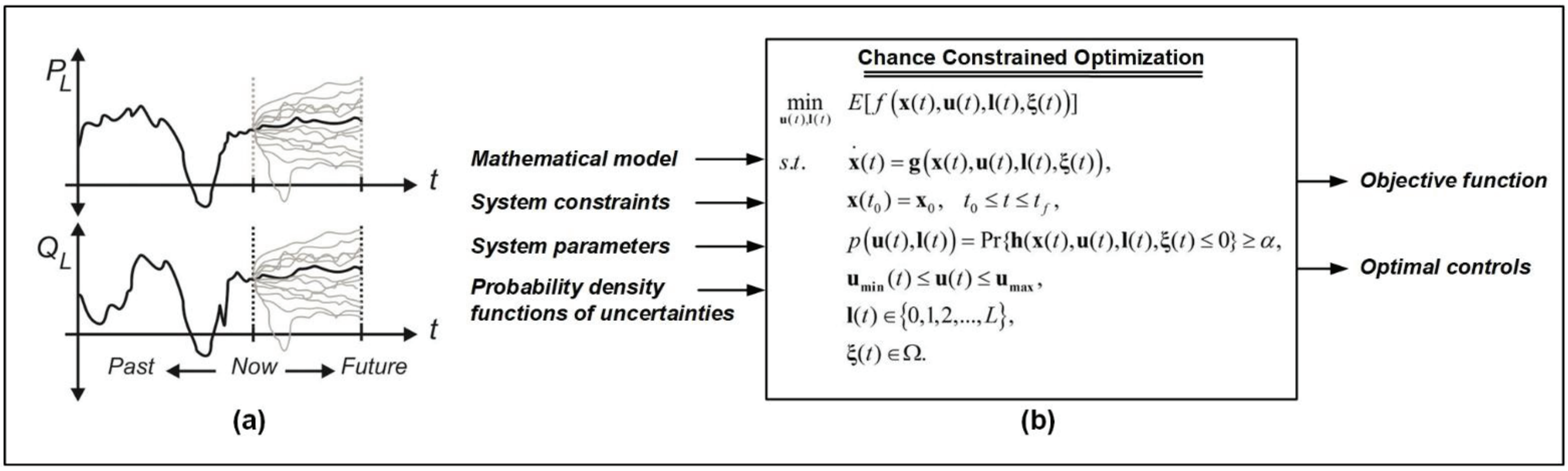

**Figure 1:** (a) Uncertain active-reactive load power of the electric machine. (b) General formulation of the chance constrained optimization problem.